\documentclass[a4paper, 11pt,twoside]{article}

\usepackage{amssymb}
\usepackage{bbm}

\newtheorem{theorem}{Theorem}
\title{Crushing candies on the line}

\begin{document}
\maketitle

\section{The candy crush game}

In this paper we investigate stability properties of a probabilistic cellular automaton. The model is based on the candy crush game (see e.g. \cite{king},\cite{wiki}).
The idea of the game is that each point (these are the candies) in
a rectangular grid has a certain color. All candies making a
horizontal or vertical monochromatic chain disappear
simultaneously and candies fall from the top to fill in the gaps.
If the resulting configuration again contains such chains, a
sequence of reactions occurs and the player is awarded bonus
points. This game inspired us to study a somewhat simplified
model, where one of the main questions is under what conditions
infinite sequences of reactions can occur.

\section{Model based on the game}

Our model is defined on $\mathbb{Z}^d$ with the Euclidean metric.
The points in $\mathbb{Z}^d$ are called sites. We choose $n$
different colors and define the color set $C =
\left\{c_1,\ldots,c_n\right\}$. Furthermore, we choose a
\emph{stability constant} $\kappa\in\mathbb{N},\kappa\geq 2$. A
\emph{coloring} or \emph{configuration} is a map
$\eta:\mathbb{Z}^d\rightarrow C$. Given the coloring, we will
define the notion of stability. Sites
$x_1,\ldots,x_m\in\mathbb{Z}^d$ satisfying
$$
x_i-x_{i-1} = x_2-x_1,\quad ||x_i-x_{i-1}||=1\quad \textrm{and}\quad \eta(x_i) = \eta(x_1)
$$
for $i = 2,\ldots,m$ are said to make a \emph{monochromatic chain}
of length $m$. Sites that are in a monochromatic chain of length
greater than or equal to $\kappa$ are called \emph{unstable}. All
other sites are stable. Given a configuration $\eta$ we define a
stability function
$\sigma_\eta:\mathbb{Z}^d\rightarrow\left\{0,1\right\}$ by
$$
\sigma_\eta(x) = \mathbbm{1}_{\left\{x \textrm{\scriptsize\ is
stable in }\eta\right\}}.
$$
A coloring $\eta$ is stable if and only if $\sigma_\eta(x) = 1$
for all $x\in\mathbb{Z}^d$. A stable configuration in which the
color of each site only depends on the parity of the sum of its
coordinates will be referred to as a chessboard coloring. The set
of all configurations will be denoted by $\Lambda$.

In our model, the dynamics are simpler than in the game. We will
use a time variable $t$, taking values in $\mathbb{N}$. The
configuration at time $t$ will be denoted by $\eta_t$, so $\eta_0$
is the initial configuration. All sites that are unstable at time
$t$ (that is, sites $x$ for which $\sigma_{\eta_t}(x)=0$) will be
recolored simultaneously and independently according to some
probability distribution $p$ on $C$ to construct $\eta_{t+1}$. For
$\eta\in\Lambda$, we will denote the random configuration that
results from recoloring the unstable sites by $R(\eta)$, so
$\eta_{t+1} = R(\eta)$.


\section{The one-dimensional case with two colors}

Let $d=1$ and $\kappa=3$. Furthermore, let the color set be coded
by $C= \left\{ 0,1 \right\}$ and choose equal recoloring
probabilities, so $p = (\frac{1}{2},\frac{1}{2})$. Then we have
the following result:
\begin{theorem}
Choose an unstable initial configuration $\eta_0$. Suppose there
exists $M\in\mathbb{N}$ such that $\sigma_{\eta_0}(x)=1$ for all
$|x|\geq M$. Then
$$
\mathbb{P}_p(\exists t_0:\eta_t \textrm{ is stable for all } t\geq t_0)=1.
$$
\end{theorem}

\textbf{Proof.} First define the number of unstable sites at time $t$:
$$
I_t = |\left\{x:\sigma_{\eta_t}(x)=0\right\}|,\quad t\geq 0.
$$
Note that $I_t$ is bounded by $2M+4t+1$ for all $t$. We will show
that $\lim_{t\rightarrow\infty}I_t = 0$ a.s. We define an upper
bounds for the probability that an unstable site is instable again
after $k$ time steps as follows:
$$
\begin{array}{lll}
p_k^I &=& \sup_{\eta\in\Lambda}\mathbb{P}(\sigma_{R^k(\eta)}(x)=0
\mid \sigma_{\eta}(x)=0),\\
p_k^{III} &=&
\sup_{\eta\in\Lambda}\mathbb{P}(\sigma_{R^k(\eta)}(x)=0 \mid
\sigma_{\eta}(x-1)=\sigma_{\eta}(x)=\sigma_{\eta}(x+1)=0).
\end{array}
$$
By translation invariance, these upper bound do not depend on $x$.
Note that $p_k^{III} \leq p_k^I$. In a similar way we define upper
bounds for the probability that a stable site is instable after
$k$ recolorings of the configuration. Since this probability
highly depends on the distance to instable regions, we condition
here on stability of a neighborhood of $x$:
$$
p_k^S(n,m) =
\sup_{\eta\in\Lambda}\mathbb{P}(\sigma_{R^k(\eta)}(x)=0 \mid
\prod_{i = -n}^{m}\sigma_{\eta}(x+i)=1,\sum_{i=-n-1}^{m+1}\sigma_{\eta}(x+i)=n+m+1),
$$
where $n,m$ are non-negative and allowed to take the value $\infty$. First we
remark that these probabilities are symmetric:
$p_k^S(n,m)=p_k^S(m,n)$. Since a stable site can not get instable
if a large enough neighborhood is stable, the probabilities
$p_k^S(n,m)$ satisfy the following property:
\begin{equation}\label{eq:pk_largenandm}
p_k^S(n,m) = 0 \quad{\rm if}\quad n,m\geq 2k.
\end{equation}
Furthermore, if there are at least $2k$ stable sites at one side,
$p_k^S(n,m)$ does not depend on the exact number of them:
\begin{equation}\label{eq:pk_forlargen}
\begin{array}{lll}
p_k^S(n_1,m) &=& p_k^S(n_2,m) = p_k^S(\infty,m) \quad{\rm if}\quad n_1,n_2\geq 2k,\\
p_k^S(n,m_1) &=& p_k^S(n,m_2) = p_k^S(n,\infty) \quad{\rm if}\quad
m_1,m_2\geq 2k.
\end{array}
\end{equation}
Let $\eta\in\Lambda$. A set $\left\{x,x+1,\ldots,x+g\right\}$ is
called a bounded stable region of size $g$ if
$\sigma_\eta(x)=\ldots=\sigma_\eta(x+g)=1$ and
$\sigma_\eta(x-1)=\sigma_\eta(x+g+1)=0$. Suppose $G$ is a bounded
stable region of size $g$, then the expected number of sites in
$G$ that get instable in $k$ steps is bounded from above by
\begin{equation}\label{eq:gapbound}
\sum_{i=1}^g p_k^S(i-1,g-i).
\end{equation}
This sum is the same for all $g>4k$, since in that case
\begin{eqnarray*}
\sum_{i=1}^g p_k^S(i-1,g-i) &=& \sum_{i=1}^{2k} p_k^S(i-1,\infty) + \sum_{i=g-2k+1}^g p_k^S(\infty,g-i)\\
&=& \sum_{i=1}^{2k} p_k^S(i-1,4k-i) + \sum_{i=2k+1}^{4k} p_k^S(i-1,4k-i)\\ &=& \sum_{i=1}^{4k} p_k^S(i-1,4k-i),
\end{eqnarray*}
by the properties (\ref{eq:pk_largenandm}) and
(\ref{eq:pk_forlargen}). If $\eta$ is such that $\sigma_\eta(x) =
1$ and $\sigma_\eta(y) = 0$ for all $y<x$, then the set
$\left\{y\in \mathbb{Z}:y<x\right\}$ is called a left-unbounded
stable region. A right-unbounded stable region is defined
similarly. For an unbounded stable region, the expected number of
sites that gets instable is at most
\begin{eqnarray*}
\sum_{i=1}^\infty p_k^S(i-1,\infty) = \sum_{i=1}^{2k}
p_k^S(i-1,\infty) = \sum_{i=2k+1}^{4k} p_k^S(4k-i,\infty).
\end{eqnarray*}
where we used (\ref{eq:pk_largenandm}) for the second equality.
Exploiting (\ref{eq:pk_forlargen}), we arrive at
\begin{eqnarray*}
\sum_{i=1}^\infty p_k^S(i-1,\infty) &=& \frac{1}{2}\sum_{i=1}^{2k}
p_k^S(i-1,4k-i) + \frac{1}{2}\sum_{i=2k+1}^{4k} p_k^S(4k-i,i-1),\\
&=& \frac{1}{2}\sum_{i=1}^{4k} p_k^S(i-1,4k-i).
\end{eqnarray*}
Since an instable region consists of at least 3 sites, at least a
third of the unstable sites has two unstable neighbors. For the
same reason, the number of bounded stable regions in $\eta_t$ is
at most $I_t/3-1$. Furthermore, since $I_t$ is finite, there are
two unbounded stable regions in $\eta_t$. Therefore
\begin{equation}\label{eq:supermartingale}
\begin{array}{lll}
\displaystyle \mathbb{E}[I_{t+k}|I_t] &\leq& \displaystyle \frac{1}{3}p_k^{III} I_t+ \frac{2}{3}p_k^{I} I_t + (\frac{I_t}{3}-1)\max_{1\leq g\leq 4k} \sum_{i=1}^g p_k^S(i-1,g-i) + \sum_{i=1}^{4k} p_k^S(i-1,4k-i)\\
&\leq& \displaystyle\left(\frac{1}{3}p_k^{III} +
\frac{2}{3}p_k^{I} +\frac{1}{3}\max_{1\leq g\leq 4k}
\sum_{i=1}^g p_k^S(i-1,g-i)\right)I_t.
\end{array}
\end{equation}
Our next goal is to show that there exists $k$ for which the
constant in front of $I_t$ is smaller than $1$. In order to do
this, we compute $p_k^I$ and $p_k^S(n,m)$ for $0\leq n,m < 2k$ for
some values of $k$.

First we take $k=1$. Let $\eta\in\Lambda$. Then
$\mathbb{P}(\sigma_{R^1(\eta)}(0)=0)$ only depends on $\eta(x)$
and $\sigma_\eta(x)$ for $-2\leq x\leq 2$. So to find the
probabilities $p_k^I$ and $p_k^S(n,m)$, we can just check all
possibilities. In the table below we listed them for the case
$\sigma_\eta(0)=0$, without loss of generality assuming that
$\eta(0)=0$ and omitting (symmetrically) equivalent cases:
$$
\begin{array}{c|c|c}
(\eta(x))_{x=-2}^{x=2} &(\sigma_\eta(x))_{x=-2}^{x=2} & \mathbb{P}(\sigma_{R^1(\eta)}(0)=0)\\
\hline
00000 & 00000 & 1/2\\
00001 & 00001 & 1/2\\
00010 & 00010 & 1/2\\
00010 & 00011 & 3/8\\
00011 & 00011 & 5/8\\
10001 & 10001 & 1/2
\end{array}
$$
It follows that $p_1^I = 5/8$ and $p_1^{III} = 1/2$. As a second
example, we compute $p_1^S(1,2)$. So here we maximize
$\mathbb{P}(\sigma_{R^1(\eta)}(0)=0)$ over configurations $\eta$
for which $\sigma_\eta(-2)=0$ and $\sigma_\eta(x)=1$ if $-1\leq
x\leq 2$. Again we assume that $\eta(0)=0$. Then the following
cases are possible:
$$
\begin{array}{c|c}
(\eta(x))_{x=-2}^{x=2} & \mathbb{P}(\sigma_{R^1(\eta)}(0)=0)\\
\hline
01001 & 0\\
01010 & 0\\
01011 & 0\\
10010 & 1/2\\
10011 & 1/2
\end{array}
$$
Therefore, $p_1^S(1,2) = 1/2$. For other values of $n$ and $m$ a
similar calculation leads to the following values of $p_1^S(n,m)$:
$$
\begin{array}{c|ccc}
&m=0&m=1&m=2\\
\hline
n=0&1/2&3/4&1/2\\
n=1&3/4&1/2&1/2\\
n=2&1/2&1/2&0
\end{array}
$$
Using (\ref{eq:pk_forlargen}) the maximal value of the sum in
(\ref{eq:supermartingale}) turns out to be equal to $2$ (for
$g=4$), whence
$$
\mathbb{E}[I_{t+1}|I_t] \leq
\left(\frac{1}{3}\cdot\frac{1}{2}+\frac{2}{3}\cdot\frac{5}{8}+\frac{1}{3}\cdot
2\right)I_t = \frac{5}{4} I_t.
$$
For other values of $k$, we can proceed analogously. The key
observation is that $\mathbb{P}(\sigma_{R^k(\eta)}(0)=0)$ is
completely determined by $\eta(x)$ and $\sigma_\eta(x)$ for
$-2k\leq x\leq 2k$. For $k=4$, we find the desired inequality
$\mathbb{E}[I_{t+k}|I_t] \leq cI_t$ with $c<1$. Below is a summary
of the computational results.
\begin{eqnarray*}
\mathbb{E}[I_{t+2}|I_t] &\leq& \textstyle\left(\frac{1}{3}\cdot\frac{29}{64}+\frac{2}{3}\cdot\frac{61}{128}+\frac{1}{3}\cdot \frac{19}{8}\right)I_t = \frac{121}{96} I_t,\\
\mathbb{E}[I_{t+3}|I_t] &\leq& \textstyle\left(\frac{1}{3}\cdot\frac{5037}{16384}+\frac{2}{3}\cdot\frac{2687}{8192}+\frac{1}{3}\cdot \frac{2495}{1024}\right)I_t = \frac{55705}{49152} I_t,\\
\mathbb{E}[I_{t+4}|I_t] &\leq& \textstyle
\left(\frac{1}{3}\cdot\frac{15371121}{67108864}+\frac{2}{3}\cdot\frac{518955}{2097152}+\frac{1}{3}\cdot
\frac{2371247}{1048576}\right)I_t = \frac{200344049}{201326592}
I_t.
\end{eqnarray*}
The values of $p_4^S(n,m)$ can be written as fractions. Their
numerators are given in the table below. Entries in the same
column have the same denominator, which is written in the second
line of the table. All denominators are powers of $2$. For example
$p_4^S(5,6)=p_4^S(6,5)=\frac{358}{2048}( =\frac{179}{1024})$.
$$
\begin{array}{c|ccccccccc}
&\scriptstyle{n=0}&\scriptstyle{n=1}&\scriptstyle{n=2}&\scriptstyle{n=3}&\scriptstyle{n=4}&\scriptstyle{n=5}&\scriptstyle{n=6}&\scriptstyle{n=7}&\scriptstyle{n=8}\\
\hline
\scriptstyle\rm{denom.}&\scriptstyle 536870912&\scriptstyle 67108864&\scriptstyle 67108864&\scriptstyle 4194304&\scriptstyle 131072&\scriptstyle 2048&\scriptstyle 512&\scriptstyle 16&\scriptstyle 1\\
&\scriptstyle (=2^{29})&\scriptstyle (=2^{26})&\scriptstyle (=2^{26})&\scriptstyle (=2^{22})&\scriptstyle (=2^{17})&\scriptstyle (=2^{11})&\scriptstyle (=2^{9})&\scriptstyle (=2^{4})&\scriptstyle (=2^{0})\\
\hline
 \scriptstyle{m=0} &\scriptstyle{ 109921252 }&&&&&&&&\\
 \scriptstyle{m=1} &\scriptstyle{ 125921345 }&\scriptstyle{  14557783 }&&&&&&&\\
 \scriptstyle{m=2} &\scriptstyle{ 116330096 }&\scriptstyle{  17164756 }&\scriptstyle{  15205719 }&&&&&&\\
 \scriptstyle{m=3} &\scriptstyle{ 131398816 }&\scriptstyle{  15993296 }&\scriptstyle{  18114414 }&\scriptstyle{   926907 }&&&&&\\
 \scriptstyle{m=4} &\scriptstyle{ 118712000 }&\scriptstyle{  17449960 }&\scriptstyle{  14234240 }&\scriptstyle{  1036746 }&\scriptstyle{   21623 }&&&&\\
 \scriptstyle{m=5} &\scriptstyle{ 128603456 }&\scriptstyle{  15641520 }&\scriptstyle{  16272448 }&\scriptstyle{   860952 }&\scriptstyle{   27214 }&\scriptstyle{   307 }&&&\\
 \scriptstyle{m=6} &\scriptstyle{ 112440704 }&\scriptstyle{  16537568 }&\scriptstyle{  12895872 }&\scriptstyle{   948720 }&\scriptstyle{   17728 }&\scriptstyle{   358 }&\scriptstyle{   47 }&&\\
 \scriptstyle{m=7} &\scriptstyle{ 119681024 }&\scriptstyle{  14745472 }&\scriptstyle{  14295552 }&\scriptstyle{   773184 }&\scriptstyle{   21120 }&\scriptstyle{   256 }&\scriptstyle{   62 }&\scriptstyle{   1 }&\\
 \scriptstyle{m=8} &\scriptstyle{ 105200384 }&\scriptstyle{  14745472 }&\scriptstyle{  11496192 }&\scriptstyle{   773184 }&\scriptstyle{   14336 }&\scriptstyle{   256 }&\scriptstyle{   32 }&\scriptstyle{   1 }&\scriptstyle{
 0}
\end{array}
$$
Now we have the inequality
$$
\mathbb{E}[I_{t+4}|I_t] \leq cI_t, \quad\quad\quad\quad{\rm with\ } c=\frac{200344049}{201326592}<1.
$$
Taking expectations,
$$
\mathbb{E}[I_{t+4}]=\mathbb{E}[\mathbb{E}[I_{t+4}|I_t]] \leq c\mathbb{E}[I_t].
$$
Therefore, for all $t\in
\mathbb{N}$, we obtain
$$
\mathbb{E}[I_{4t}] \leq c^{t} \mathbb{E}[I_0] \leq c^t(2M+1).
$$
Since $I_{4t}$ is positive integer-valued and the events $\left\{I_{4t}\geq 1\right\}$ are decreasing, we get
\begin{eqnarray*}
\mathbb{P}(\exists t:I_{4t}=0) &=& \mathbb{P}\left(\bigcup_{t=0}^\infty \left\{I_{4t}=0\right\}\right)
= 1-\mathbb{P}\left(\bigcap_{t=0}^\infty \left\{I_{4t}\geq 1\right\}\right)\\
&=& 1-\lim_{t\rightarrow\infty} \mathbb{P}\left(\left\{I_{4t}\geq 1\right\}\right)
\geq 1-\lim_{t\rightarrow\infty} \mathbb{E}[I_{4t}]\\
&\geq& 1-\lim_{t\rightarrow\infty} c^t(2M+1) = 1,
\end{eqnarray*}
where we used Markov's inequality. If $I_{4t_0}=0$, then $\eta_{t}$ is stable for all $t\geq 4t_0$. Therefore
$$
\mathbb{P}(\exists t_0:\eta_t \textrm{ is stable for all } t\geq t_0)=1.
$$
\hfill$\Box$

\end{document}